\documentclass[12pt]{article}
\usepackage{amsmath,amsthm, amssymb}
\usepackage{amssymb,latexsym}

\usepackage{eso-pic}
\usepackage{lipsum}

\newtheorem{theorem}{Theorem}

\newtheorem{lemma}{Lemma}

\begin{document}
\baselineskip=17pt

\AddToShipoutPictureBG*{
  \AtPageUpperLeft{
  \hspace{\paperwidth}
   \raisebox{-4\baselineskip}{
   \makebox[-100pt][r]{
     This is the author's version of the paper.  The final publication has appeared in}

     \raisebox{-1.2\baselineskip}{
      \makebox[-166pt][r]{
     Proc. Techn. Univ.-Sofia, \textbf{66}, 3, (2016), 55 -- 64.}
}}}}

\title{\bf Arithmetic progressions of three prime numbers
with two primes of the form $\mathbf{p=x^2+y^2+1}$}

\author{\bf S. I. Dimitrov}

\date{\bf2016}
\maketitle

\begin{abstract}
In the present  paper we prove that there exist infinitely many arithmetic progressions of three different primes
$p_1,p_2,p_3=2p_2-p_1$ such that $p_1=x_1^2 + y_1^2 +1$, $p_2=x_2^2 + y_2^2 +1$.
\medskip

\end{abstract}

\textbf{Notations.} The letter $p$, with or without subscript, will always denote prime numbers.
By $\varepsilon$ we denote an arbitrary small positive number, not the same in all appearances.
$X$ is an sufficiently large positive number and $\mathcal{L}=\log X$. We denote by $\mathcal{J}$
the set of all subintervals of the interval $(X/2,X]$ and if $J_1,J_2\in \mathcal{J}$ then
$\textbf{J}=\langle J_1,J_2\rangle$ is the corresponding ordered pair. Respectively
$\textbf{d}=\langle d_1,d_2\rangle$ is an two-dimensional vector with integer components $d_1,d_2$
and in particular $\textbf{1}=\langle1,1\rangle$. We denote by $(m,n)$ the greatest common divisor
of $m$ and $n$. As usual $\varphi(d)$ is Euler's function; $\tau(d)$ is the number of positive
divisors of $d$; $r(d)$ is the number of solutions of the equation $d=m_1^2+m_2^2$ in integers
$m_j$; $\chi(d)$ is the non-principal character modulo 4 and $L(s,\chi)$ is the corresponding
Dirichlet's $L$-function.
\section{Introduction and statement of the result.}
\indent

In 1939 Van der Corput \cite{Corput} proved that there exist infinitely many arithmetic progressions of three different primes.
On the other hand Linnik \cite{Linnik} has proved that there exist
infinitely many prime numbers of the form $p=x^2 + y^2 +1$, where $x$ and $y$ -- integers.
More precisely he has proved the asymptotic formula
\begin{equation}\label{Linnik}
\sum_{p\leq X}r(p-1)=\pi X\mathcal{L}^{-1}\prod_{p>2}\bigg(1+\frac{\chi(p)}{p(p-1)}\bigg)+
\mathcal{O}\left(X\mathcal{L}^{-1-\theta_0}(\log\mathcal{L})^7\right),
\end{equation}
where
\begin{equation}\label{theta0}
\theta_0=\frac{1}{2}-\frac{1}{4}e\log2=0.0289...
\end{equation}
We couple these two theorems by proving the following theorem.

Define
\begin{equation} \label{R}
R(X)=\sum\limits_{X/2<p_1,p_2,p_3\leq X\atop{p_1+p_3=2p_2}}r(p_1-1)r(p_2-1)\log p_1\log p_2\log p_3
\end{equation}
and
\begin{equation}\label{sigma0}
\sigma_0=2\prod\limits_{p>2}\left(1-\frac{1}{(p-1)^2}\right),
\end{equation}
\begin{equation}\label{SigmaR}
\mathfrak{S}_R=\pi^2\sigma_0\prod\limits_p\left(1+\chi(p)\frac{2p^2+p\chi(p)-4p+2\chi(p)}{p^2(p-1)(p-2)}\right).
\end{equation}
\begin{theorem}We have the following asymptotic formula
\begin{equation}\label{RAs}
R(X)=\frac{1}{8}\mathfrak{S}_RX^2+\mathcal{O}\big(X^2\mathcal{L}^{-\theta_0}(\log\mathcal{L})^7\big)\,,
\end{equation}
where $\theta_0$ is denoted by  \eqref{theta0}.
\end{theorem}
We use sieve methods to pick out primes of the form $x^2 + y^2 +1$ and the circle method to pick out primes satisfying
$p_1+p_3=2p_2$.
Recently Matom\"{a}ki  \cite{Mato} and Tolev \cite{Tolev1}  have obtained a similar results related to a binary Goldbach problem.
Tolev \cite{Tolev2} has also proved that every sufficiently large odd integer can be represented as a sum
of three primes, two of which of the form $x^2 + y^2 +1$.
Our argument is a modification of Tolev's \cite{Tolev2} argument.

\section{Some lemmas.}
\indent

First we consider the binary Goldbach problem with one prime variable lying in a given interval
and belonging to an arithmetic progression.
Suppose that $n\leq 2X$, let $d$ and $l$ be integers with $(d,l)=1$ and let
$J\in\mathcal{J}$.

Denote
\begin{align}\label{lambda}
&\lambda(n)=\prod\limits_{p|n\atop{p>2}}\frac{p-1}{p-2}\,;\\
\label{Sigma}
&\mathfrak{S}_{d,l}(n)=\begin{cases}\sigma_0\lambda(nd)\;\;\;\;\mbox{if}\;\;(d,n-l)=1\;\;\mbox{and}\;\; 2\mid n,\\
 0\;\;\;\; \;\;\;\;\;\;\;\;\;\;\mbox{otherwise}\,;\\
 \end{cases}\\
\label{I1}
&I^{(1)}_{d,l}(n,X,J)=\sum\limits_{X/2<p_1\leq X\atop{2p_2-p_1=n\atop{p_2\equiv l\;(\textmd{mod}\,d)
\atop{p_2\in J}}}}\log p_1\log p_2\,;\\
\label{Phi1}
&\Phi^{(1)}(n,X,J)=\sum\limits_{X/2<m_1\leq X\atop{2m_2-m_1=n\atop{m_2\in J}}}1\,;\\
\label{Delta1}
&\Delta_{d,l}^{(1)}(n,X,J)=I^{(1)}_{d,l}(n,X,J)-\frac{\mathfrak{S}_{d,l}(n)}{\varphi(d)}\Phi^{(1)}(n,X,J)\,;\\
\label{I2}
&I^{(2)}_{d,l}(n,X,J)=\sum\limits_{X/2<p_1\leq X\atop{p_1+p_2=n\atop{p_2\equiv l\;(\textmd{mod}\,d)
\atop{p_2\in J}}}}\log p_1\log p_2\,;\\
\label{Phi2}
&\Phi^{(2)}(n,X,J)=\sum\limits_{X/2<m_1\leq X\atop{m_1+m_2=n\atop{m_2\in J}}}1\,;\\
\label{Delta2}
&\Delta_{d,l}^{(2)}(n,X,J)=I^{(2)}_{d,l}(n,X,J)-\frac{\mathfrak{S}_{d,l}(n)}{\varphi(d)}\Phi^{(2)}(n,X,J)\,.
\end{align}
If $J=(X/2,X]$ then we write for simplicity $I^{(j)}_{d,l}(n,X)$,\;$\Phi^{(j)}(n,X)$ and
$\Delta_{d,l}^{(j)}(n,X)$ for $j=1,2$.

The following Bombieri-Vinogradov type result gives the arithmetical information needed for
the applications of the sieve.
\begin{lemma}\label{Laporta}
For any constant $A>0$ there exists $B=B(A)>0$ such that
\begin{equation*}
\sum\limits_{d\le\sqrt{X}\mathcal{L}^{-B}}\max\limits_{(d,l)=1}\max\limits_{J\in\mathcal{J}}
\sum\limits_{n\le X}\big|\Delta_{d,l}^{(j)}(n,X,J)\big|\ll X^2\mathcal{L}^{-A},\;j=1,2.
\end{equation*}
\end{lemma}
This lemma is a very similar to the result of Laporta \cite{Laporta}. This author studies the equation $p_1-p_2=n$
and without the condition $p_1\in J$. However, inspecting the arguments presented
in \cite{Laporta}, the reader will readily see that the proof of Lemma \ref{Laporta} can be obtained in the same way.

Next we consider the equation $p_1+p_3=2p_2$  with two primes from arithmetic progressions and belonging to given
intervals. Suppose that $\textbf{d}=\langle d_1,d_2\rangle$ and  $\textbf{l}=\langle l_1,l_2\rangle$
are two-dimensional vectors with integer components such that $(d_i,l_i)=1,\;i=1,2$ and let
$\textbf{J}=\langle J_1,J_2\rangle$ be a pair of intervals $J_1,J_2\in \mathcal{J}$.

Denote
\begin{align}\label{I3}
&I^{(3)}_{\textbf{d},\textbf{l}}(X,\textbf{J})=\sum\limits_{X/2<p_3\leq X\atop{p_1+p_3=2p_2\atop{p_i\equiv l_i\,(\textmd{mod}\,d_i)
\atop{p_i\in J_i\,,i=1,2}}}}\log p_1\log p_2\log p_3\,;\\
\label{Phi3}
&\Phi^{(3)}(X,\textbf{J})=\sum\limits_{X/2<m_3\leq X\atop{m_1+m_3=2m_2\atop{m_i\in J_i\,,i=1,2}}}1\,.
\end{align}
Next we define $\mathfrak{S}^{(3)}_{\textbf{d},\textbf{l}}$ in the following way. Consider the sets of primes
\begin{align*}
&\mathcal{A}=\{p:\;p\nmid d_1d_2\}\,;\\
&\mathcal{B}=\{p:\;p\mid d_1,\;p\nmid d_2\}\cup\{p:\;p\nmid d_1,\;p\mid d_2\}\,;\\
&\mathcal{C}=\{p:\;p\mid d_1,\;p\mid d_2,\;p\mid(l_1-2l_2)\}\,;\\
&\mathcal{D}=\{p:\;p\mid d_1,\;p\mid d_2,\;p\nmid(l_1-2l_2)\}\,.
\end{align*}
If $\mathcal{C}\neq\emptyset$ then we assume that
\begin{equation}\label{zero}
\mathfrak{S}^{(3)}_{\textbf{d},\textbf{l}}=0\,.
\end{equation}
If $\mathcal{C}=\emptyset$ then we put
\begin{equation}\label{Sigma3}
\mathfrak{S}^{(3)}_{\textbf{d},\textbf{l}}=2\prod\limits_{p\in\mathcal{A}\cup\mathcal{B}}
\left(1-\frac{1}{(p-1)^2}\right)\prod\limits_{p\in\mathcal{D}}\left(1+\frac{1}{p-1}\right)\,.
\end{equation}
We also define
\begin{equation}\label{Delta3}
\Delta_{\textbf{d},\textbf{l}}^{(3)}(X,\textbf{J})=I^{(3)}_{\textbf{d},\textbf{l}}(X,\textbf{J})
-\frac{\mathfrak{S}^{(3)}_{\textbf{d},\textbf{l}}}
{\varphi(d_1)\varphi(d_2)}\Phi^{(3)}(X,\textbf{J})\,.
\end{equation}
If $J_1=J_2=(X/2,X]$ then we write for simplicity $I^{(3)}_{\textbf{d},\textbf{l}}(X)$,\;$\Phi^{(3)}(X)$
and $\Delta_{\textbf{d},\textbf{l}}^{(3)}(X)$.

The next lemma is analogous to Lemma \ref{Laporta} and also gives the arithmetical information needed for
the applications of the sieve.
\begin{lemma}\label{Peneva,Tolev}
For any constant $A>0$ there exists $B=B(A)>0$ such that
\begin{equation*}
\sum\limits_{d_1\le\sqrt{X}\mathcal{L}^{-C}}\sum\limits_{d_2\le\sqrt{X}\mathcal{L}^{-C}}
\max\limits_{(d_i,l_i)=1\atop{i=1,2}}\max\limits_{J_i\in\mathcal{J}\atop{i=1,2}}
\big|\Delta_{\textbf{d},\textbf{l}}^{(3)}(X,\textbf{J})\big|\ll X^2\mathcal{L}^{-A}\,.
\end{equation*}
\end{lemma}

This assertion is a slight generalization to the theorem of Peneva and Tolev \cite{PenevaTolev}.
In their theorem $(d_1d_2,2)=1$, $l_1=l_2=-2$ and there are no conditions $p_i\in J_i$.
But it is not hard to verify that the method used in this paper implies also the correctness of
Lemma \ref{Peneva,Tolev}.

Several times we shall use the following
\begin{lemma}\label{j}
Suppose that $j\in\{-1,1\}$ and let $d,l,m$ be natural numbers. Then the quantities $\mathfrak{S}_{4m,1+jm}(n)$
and $\mathfrak{S}^{(3)}_{\langle d,4m\rangle,\langle l,1+jm\rangle}$ do not depend on $j$.
\end{lemma}
The proof is a immediate consequence from the definitions of $\mathfrak{S}_{d,l}(n)$ and
 $\mathfrak{S}^{(3)}_{\textbf{d},\textbf{l}}$.

The next two lemmas are due to C. Hooley.
\begin{lemma}\label{Hooley1}
For any constant $\omega>0$ we have
\begin{equation*}
\sum\limits_{p\leq X}\bigg|\sum\limits_{d|p-1
\atop{\sqrt{X}\mathcal{L}^{-\omega}<d<\sqrt{X}\mathcal{L}^{\omega}}}\chi(d)\bigg|
\ll X\mathcal{L}^{-1-\theta_0}(\log\mathcal{L})^5\,,
\end{equation*}
where $\theta_0$ is defined by \eqref{theta0}. The constant in the Vinogradov symbol depends only on $\omega>0$.
\end{lemma}
\begin{lemma}\label{Hooley2}
For any constant $\omega>0$ we have
\begin{equation*}
\sum\limits_{p\leq X}\bigg|\sum\limits_{d|p-1
\atop{\sqrt{X}\mathcal{L}^{-\omega}<d<\sqrt{X}\mathcal{L}^{\omega}}}\chi(d)\bigg|^2
\ll X\mathcal{L}^{-1}(\log\mathcal{L})^7\,,
\end{equation*}
where the constant in the Vinogradov symbol depends on $\omega>0$.
\end{lemma}
The proofs of very similar results are available in (\cite{Hooley}, Ch.5).

The next lemma is due to Tolev (\cite{Tolev2}, Lemma 6).
\begin{lemma}\label{Tolev}Let $n$ be an integer satisfying $1\leq n\leq X$. Suppose that $\omega>0$ is a constant
and let $\mathcal{P}=\mathcal{P}_\omega(X)$ be the set of primes $p\leq X$ such that $p-1$ has a divisor lying between $\sqrt{X}\mathcal{L}^{-\omega}<d<\sqrt{X}\mathcal{L}^{\omega}$. Then we have
\begin{equation*}
\sum\limits_{p_1+p_2=n\atop{p_1\in \mathcal{P}}}1\ll X\mathcal{L}^{-2-2\theta_0}
(\log\mathcal{L})^6\,,
\end{equation*}
where $\theta_0$ is defined by \eqref{theta0}. The constant in the Vinogradov's symbol depends only on $\omega>0$.
\end{lemma}
Arguing as in Lemma 6 we obtain similar result.
\begin{lemma}\label{Tol-Dim}Let $n$ be an integer satisfying $1\leq n\leq X$. Suppose that $\omega>0$ is a constant
and let $\mathcal{P}=\mathcal{P}_\omega(X)$ be the set of primes $p\leq X$ such that $p-1$ has a divisor lying between $\sqrt{X}\mathcal{L}^{-\omega}<d<\sqrt{X}\mathcal{L}^{\omega}$. Then we have
\begin{equation*}
\sum\limits_{2p_1-p_2=n\atop{p_1\in \mathcal{P}}}1\ll X\mathcal{L}^{-2-2\theta_0}
(\log\mathcal{L})^6\,,
\end{equation*}
where $\theta_0$ is defined by \eqref{theta0}. The constant in the Vinogradov's symbol depends only on $\omega>0$.
\end{lemma}
\section{Proof of Theorem 1.}
\indent

\textbf{Beginning.}
Denote
\begin{equation}\label{D}
D=\sqrt{X}\mathcal{L}^{-B(10)-C(10)}\,,
\end{equation}
where $B(A)$ and $C(A)$ are specified respectively in Lemma \ref{Laporta} and Lemma \ref{Peneva,Tolev}.\\
Obviously
\begin{equation}\label{rm}
r(m)=4\sum\limits_{d\mid m}\chi(d)= 4\big(r_1(m)+r_2(m)+r_3(m)\big)\,,
\end{equation}
where
\begin{equation}\label{r1r2r3}
r_1(m)=\sum\limits_{d\mid m\atop{d\leq D}}\chi(d)\,,\;\;\;\;\;
r_2(m)=\sum\limits_{d\mid m\atop{D<d<X/D}}\chi(d)\,,\;\;\;\;\;
r_3(m)=\sum\limits_{d\mid m\atop{d\geq X/D}}\chi(d)\,.
\end{equation}
Using (\ref{R}) and (\ref{rm}) we find
\begin{equation}\label{R16}
R(X)=16\sum\limits_{1\leq i,j\leq3}R_{i,j}(X)\,,
\end{equation}
where
\begin{equation} \label{Rij}
R_{i,j}(X)=\sum\limits_{X/2<p_1,p_2,p_3\leq X\atop{p_1+p_3=2p_2}}r_i(p_1-1)r_j(p_2-1)\log p_1\log p_2\log p_3\,.
\end{equation}
We shall prove that the main term in (\ref{RAs}) comes from $R_{1,1}(X)$ and the other sums $R_{i,j}(X)$
contribute only to the remainder term.

\vspace{1mm}
\noindent\textbf{The estimation of} $\textbf{R}_{\textbf{1},\textbf{1}}\textbf{(X).}$
Using (\ref{Phi3}), (\ref{Delta3}), (\ref{r1r2r3}) and (\ref{Rij}) we obtain
\begin{equation}\label{R11}
R_{1,1}(X)=\sum\limits_{d_1,d_2\leq D}\chi(d_1)\chi(d_2)I^{(3)}_{\textbf{d},\textbf{1}}(X)
=R^\prime_{1,1}(X)+R^*_{1,1}(X)\,,
\end{equation}
where
\begin{align}\label{R'}
&R^\prime_{1,1}(X)=\Phi^{(3)}(X)\sum\limits_{d_1,d_2\leq D}\frac{\chi(d_1)\chi(d_2)}
{\varphi(d_1)\varphi(d_2)}\mathfrak{S}^{(3)}_{\textbf{d},\textbf{1}}\,,\\
\label{R11*}
&R^*_{1,1}(X)=\sum\limits_{d_1,d_2\leq D}\chi(d_1)\chi(d_2)\Delta_{\textbf{d},\textbf{1}}^{(3)}(X)\,.
\end{align}
From (\ref{D}), (\ref{R11*}) and Lemma \ref{Peneva,Tolev} we have that
\begin{equation}\label{R11*est}
R^*_{1,1}(X)\ll X^2\mathcal{L}^{-1}\,.
\end{equation}
Consider $R^\prime_{1,1}(X)$. From (\ref{Phi3}) and (\ref{R'}) it follows that
\begin{equation}\label{R11'Gamma}
R^\prime_{1,1}(X)=\frac{1}{8}X^2\Gamma(X)+\mathcal{O}\left(X^{1+\varepsilon}\right)\,,
\end{equation}
where
\begin{equation}\label{Gamma}
\Gamma(X)=\sum\limits_{d_1,d_2\leq D}\frac{\chi(d_1)\chi(d_2)}
{\varphi(d_1)\varphi(d_2)}\mathfrak{S}^{(3)}_{\textbf{d},\textbf{1}}\,.
\end{equation}
We shall find an asymptotic formula for $\Gamma(X)$.\\
Using (\ref{zero}), (\ref{Sigma3}) and (\ref{Gamma}) we get
\begin{equation}\label{Gammaf}
\Gamma(X)=\sigma_0\sum\limits_{d\leq D}\frac{\chi(d)}{\varphi(d)}\sum\limits_{t\leq D}f_d(t)\,,
\end{equation}
where
\begin{equation}\label{fdt}
f_d(t)=\frac{\chi(t)}{\varphi(t)}\prod\limits_{p\mid(d,t)}\frac{p-1}{p-2}\,.
\end{equation}
First we estimate the sum over $t$ in (\ref{Gammaf}). From (\ref{fdt}) we have that
\begin{equation}\label{fdtest}
f_d(t)\ll(\log\mathcal{L})t^{-1}\log\log(10t)
\end{equation}
with absolute constant in the Vinogradov's symbol. Thus the corresponding Dirichlet series
\begin{equation*}
F_d(s)=\sum\limits_{t=1}^\infty f_d(t)t^{-s}
\end{equation*}
is absolutely convergent in $Re(s)>0$. On the other hand $f_d(t)$ is a multiplicative with
respect to $t$ and applying Euler's identity we find
\begin{equation}\label{FT}
F_d(s)=\prod\limits_pT_d(p,s)\,,\;\;\;\;T_d(p,s)=1+\sum\limits_{l=1}^\infty f_d(p^l)p^{-ls}\,.
\end{equation}
From (\ref{fdt}) and (\ref{FT}) it follows that
\begin{equation*}
T_d(p,s)=\left(1-\chi(p)p^{-s-1}\right)^{-1}\left(1+\chi(p)p^{-s-1}Y_d(p)\right)\,,
\end{equation*}
where
\begin{equation}\label{Ydp}
Y_d(p)=\begin{cases}(p-1)^{-1}\;\;\;\;\;\mbox{if}\;\; p\nmid d,\\
2(p-2)^{-1}\;\;\;\mbox{if}\;\; p\mid d\,.
\end{cases}
\end{equation}
Hence we obtain
\begin{equation}\label{Fds}
F_d(s)=L(s+1,\chi)\prod\limits_p \left(1+\chi(p)p^{-s-1}Y_d(p)\right)\,.
\end{equation}
From this formula it follows that $F_d(s)$ has an analytic continuation to $Re(s)>-1$.
Using (\ref{Ydp}), (\ref{Fds}) and the simplest bound for $L(s+1,\chi)$ we find
\begin{equation}\label{Fdsest}
F_d(s)\ll X^{1/6}\;\;\;\mbox{for}\;\;Re(s)\geq-\frac{1}{2}\,,\;\;|Im(s)|\leq X\,.
\end{equation}
We apply Perron's formula given at Tenenbaum (\cite{Tenenbaum}, Chapter II.2) and also (\ref{fdtest}) to
get
\begin{equation}\label{sumfdt}
\sum\limits_{t\leq D}f_d(t)=\frac{1}{2\pi\imath}\int\limits_{\kappa-\imath X}^{\kappa+\imath X}F_d(s)\frac{D^s}{s}ds
+\mathcal{O}\left(\sum\limits_{t=1}^\infty\frac{X^\varepsilon D^\kappa\log\log(10t)}{t^{1+\kappa}
\left(1+X\left|\log\frac{D}{t}\right|\right)}\right)\,.
\end{equation}
where $\kappa=1/10$. It is easy to see that the error term above is $\mathcal{O}\left(X^{-1/20}\right)$.
Applying the residue theorem we see that the main term is equal to
\begin{equation*}
F_d(0)+\frac{1}{2\pi\imath}\left(\int\limits_{1/10-\imath X}^{-1/2-\imath X}+
\int\limits_{-1/2-\imath X}^{-1/2+\imath X}+\int\limits_{-1/2+\imath X}^{1/10+\imath X}\right)F_d(s)\frac{D^s}{s}ds\,.
\end{equation*}
From (\ref{Fdsest}) it follows that the contribution from the above integrals is $\mathcal{O}\left(X^{-1/20}\right)$.\\
Hence
\begin{equation}\label{sumfdtest}
\sum\limits_{t\leq D}f_d(t)=F_d(0)+\mathcal{O}\left(X^{-1/20}\right)\,.
\end{equation}
Using (\ref{Fds}) we obtain
\begin{equation}\label{Fd0}
F_d(0)=\frac{\pi}{4}\prod\limits_p\left(1+\frac{\chi(p)}{p}Y_d(p)\right)\,.
\end{equation}
Having in mind (\ref{Gammaf}), (\ref{Ydp}), (\ref{sumfdtest}) and (\ref{Fd0}) we establish that
\begin{equation}\label{Gammag}
\Gamma(X)=\frac{\pi}{4}\sigma_0\mathfrak{X}\sum\limits_{d\leq D}g(d)+\mathcal{O}\left(X^{\varepsilon-1/20}\right)\,.
\end{equation}
where $\sigma_0$ is defined by (\ref{sigma0}),
\begin{equation}\label{Xi}
\mathfrak{X}=\prod\limits_p\left(1+\frac{\chi(p)}{p(p-1)}\right)\,,
\end{equation}
and
\begin{equation}\label{gd}
g(d)=\frac{\chi(d)}{\varphi(d)}\prod\limits_{p\mid d}\frac{1+\frac{2\chi(p)}{p(p-2)}}
{1+\frac{\chi(p)}{p(p-1)}}\,.
\end{equation}
Obviously $g(d)$ is multiplicative with respect to $d$ and satisfies
\begin{equation}\label{gdest}
g(d)\ll d^{-1}\log\log(10d),
\end{equation}
where the constant in Vinogradov's symbol is absolute. Thus the Dirichlet series
\begin{equation*}
G(s)=\sum\limits_{d=1}^\infty g(d)d^{-s}
\end{equation*}
is absolutely convergent in $Re(s)>0$ and applying the Euler's identity we find
\begin{equation}\label{GH}
G(s)=\prod\limits_pH(p,s)\,,\;\;\;\;H(p,s)=1+\sum\limits_{l=1}^\infty g(p^l)p^{-ls}\,.
\end{equation}
From (\ref{gd}) and (\ref{GH}) we get
\begin{equation*}
H(p,s)=\left(1-\chi(p)p^{-s-1}\right)^{-1}\left(1+\chi(p)p^{-s-1}K(p)\right)\,,
\end{equation*}
where
\begin{equation}\label{Kp}
K(p)=\frac{p^2+p\chi(p)-2p+2\chi(p)}{p^3-3p^2+p\chi(p)+2p-2\chi(p)}\,.
\end{equation}
This implies
\begin{equation}\label{Gs}
G(s)=L(s+1,\chi)\prod\limits_p \left(1+\chi(p)p^{-s-1}K(p)\right)\,.
\end{equation}
We see that $G(s)$ has an analytic continuation to $Re(s)>-1$ and
\begin{equation}\label{Gsest}
G(s)\ll X^{1/6}\;\;\;\mbox{for}\;\;Re(s)\geq-\frac{1}{2}\,,\;\;|Im(s)|\leq X\,.
\end{equation}
Applying Perron's formula and proceeding as above we obtain
\begin{equation}\label{sumgdest}
\sum\limits_{d\leq D}g(d)=G(0)+\mathcal{O}\left(X^{-1/20}\right)=
\frac{\pi}{4}\prod\limits_p\left(1+\frac{\chi(p)}{p}K(p)\right)+\mathcal{O}\left(X^{-1/20}\right)\,.
\end{equation}
Using (\ref{Gammag}), (\ref{Xi}), (\ref{Kp}) and (\ref{sumgdest}) we get
\begin{equation}\label{Gammaest}
\Gamma(X)=\frac{1}{16}\mathfrak{S}_R+\mathcal{O}\left(X^{\varepsilon-1/20}\right)\,.
\end{equation}
where $\mathfrak{S}_R$ is defined by (\ref{SigmaR}).

Now bearing in mind (\ref{R11}), (\ref{R11*est}), (\ref{R11'Gamma}) and (\ref{Gammaest}) we find
\begin{equation}\label{R11est}
R_{1,1}(X)=\frac{1}{128}\mathfrak{S}_RX^2+\mathcal{O}\left(X^2\mathcal{L}^{-1}\right)\,.
\end{equation}
\textbf{The estimation of} $\textbf{R}_{\textbf{1},\textbf{2}}\textbf{(X).}$
Using (\ref{I2}) -- (\ref{Delta2}), (\ref{r1r2r3}) and (\ref{Rij}) we obtain
\begin{equation}\label{R12}
R_{1,2}(X)=\sum\limits_{X/2<p\leq X}r_2(p-1)\log p\sum\limits_{d\leq D}\chi(d)I^{(2)}_{d,1}(2p,X)=
R^\prime_{1,2}(X)+R^*_{1,2}(X)\,,
\end{equation}
where
\begin{align}\label{R12'}
&R^\prime_{1,2}(X)=\sum\limits_{X/2<p\leq X}r_2(p-1)\log p\sum\limits_{d\leq D}
\frac{\chi(d)}{\varphi(d)}\mathfrak{S}_{d,1}(2p)\Phi^{(2)}(2p,X)\,,\\
\label{R12*}
&R^*_{1,2}(X)=\sum\limits_{X/2<p\leq X}r_2(p-1)\log p\sum\limits_{d\leq D}\chi(d)\Delta^{(2)}_{d,1}(2p,X)\,.
\end{align}
From (\ref{r1r2r3}), (\ref{R12*}) and Cauchy's inequality it follows
\begin{align}\label{R12*est1}
R^*_{1,2}(X)&\ll\mathcal{L}\sum\limits_{X/2<p\leq X}\tau(p-1)
\sum\limits_{d\leq D}\big|\Delta^{(2)}_{d,1}(2p,X)\big|\nonumber\\
&\ll\mathcal{L}\sum\limits_{X<n\leq2X}\tau(n)\sum\limits_{d\leq D}\big|\Delta^{(2)}_{d,1}(n,X)\big|\nonumber\\
&\ll\mathcal{L}\left(\sum\limits_{X<n\leq2X}\sum\limits_{d\leq D}\tau^2(n)\big|\Delta^{(2)}_{d,1}(n,X)\big|\right)^{1/2}
\left(\sum\limits_{X<n\leq2X}\sum\limits_{d\leq D}\big|\Delta^{(2)}_{d,1}(n,X)\big|\right)^{1/2}\nonumber\\
&=\mathcal{L}U^{1/2}V^{1/2}\,,
\end{align}
say. We use the trivial estimation $\Delta^{(2)}_{d,1}(n,X)\ll\mathcal{L}^2Xd^{-1}$ and the well-known inequality
$\sum_{n\leq y}\tau^2(n)\ll y\log^3y$ to get
\begin{equation}\label{Uest}
U\ll X^2\mathcal{L}^6\,.
\end{equation}
In order to estimate $V$ we use (\ref{D}) and Lemma \ref{Laporta} to obtain
\begin{equation}\label{Vest}
V\ll X^2\mathcal{L}^{-10}\,.
\end{equation}
From (\ref{R12*est1}) -- (\ref{Vest})  it follows that
\begin{equation}\label{R12*est2}
R^*_{1,2}(X)\ll X^2\mathcal{L}^{-1}\,.
\end{equation}

Consider now $R^\prime_{1,2}(X)$. Having in mind (\ref{lambda}), (\ref{Sigma}) and (\ref{R12'}) we find
\begin{equation}\label{R12'est1}
R^\prime_{1,2}(X)=\sigma_0\sum\limits_{X/2<p\leq X}r_2(p-1)\Phi^{(2)}(2p,X)\lambda(2p)\log p
\sum\limits_{d\leq D\atop{(d,2p-1)=1}}\frac{\chi(d)\lambda(d)}{\varphi(d)\lambda((d,2p))}\,.
\end{equation}
For the sum over $d$ according to (\cite{Tolev1}, Section 3.2) we have the bound $\sum_d\ll\log\mathcal{L}$.
Therefore, using also (\ref{sigma0}) and (\ref{Phi2}) we get
\begin{equation}\label{R12'est2}
R^\prime_{1,2}(X)\ll X\mathcal{L}(\log\mathcal{L})^2\sum\limits_{X/2<p\leq X}|r_2(p-1)|\,.
\end{equation}
Bearing in mind (\ref{D}), (\ref{r1r2r3}), (\ref{R12'est2}) and Lemma \ref{Hooley1} we obtain
\begin{equation}\label{R12'est3}
R^\prime_{1,2}(X)\ll X^2\mathcal{L}^{-\theta_0}(\log\mathcal{L})^7\,.
\end{equation}

Finally from  (\ref{R12}), (\ref{R12*est2}) and  (\ref{R12'est3}) we find
\begin{equation}\label{R12est}
R_{1,2}(X)\ll X^2\mathcal{L}^{-\theta_0}(\log\mathcal{L})^7\,.
\end{equation}
\textbf{The estimation of} $\textbf{R}_{\textbf{1},\textbf{3}}\textbf{(X).}$
From (\ref{I3}), (\ref{r1r2r3}) and (\ref{Rij}) we have
\begin{align*}
R_{1,3}(X)&=\sum\limits_{X/2<p_1,p_2,p_3\leq X\atop{p_1+p_3=2p_2}}\log p_1\log p_2\log p_3
\sum\limits_{d\mid p_1-1\atop{d\leq D}}\chi(d)\sum\limits_{m\mid p_2-1\atop{\frac{p_2-1}{m}\geq X/D}}
\chi\left(\frac{p_2-1}{m}\right)\\
&=\sum\limits_{d\leq D\atop{m<D\atop{2\mid m}}}\chi(d)\sum\limits_{j=\pm1}\chi(j)
I^{(3)}_{\langle d,4m\rangle,\langle1,1+jm\rangle}(X,\langle(X/2,X],J_m\rangle)\,,
\end{align*}
where $J_m=\big(\max\{1+mX/D,X/2\},X\big]$. From (\ref{Delta3}) we find
\begin{equation}\label{R13}
R_{1,3}(X)=R^\prime_{1,3}(X)+R^*_{1,3}(X)\,,
\end{equation}
where
\begin{align}\label{R13'}
&R^\prime_{1,3}(X)=\sum\limits_{d\leq D\atop{m<D\atop{2\mid m}}}\frac{\chi(d)\Phi^{(3)}(X,\langle(X/2,X],J_m\rangle)}{\varphi(d)\varphi(4m)}\sum\limits_{j=\pm1}\chi(j)
\mathfrak{S}^{(3)}_{\langle d,4m\rangle,\langle1,1+jm\rangle}\,,\\
\label{R13*}
&R^*_{1,3}(X)=\sum\limits_{d\leq D\atop{m<D\atop{2\mid m}}}\chi(d)\sum\limits_{j=\pm1}\chi(j)
\Delta^{(3)}_{\langle d,4m\rangle,\langle1,1+jm\rangle}(X,\langle(X/2,X],J_m\rangle)\,.
\end{align}
From (\ref{D}) and Lemma \ref{Peneva,Tolev} we obtain
\begin{equation}\label{R13*est}
R^*_{1,3}(X)\ll X^{2}\mathcal{L}^{-1}\,.
\end{equation}

Consider $R^\prime_{1,3}(X)$. According to Lemma \ref{j} the expression
$\mathfrak{S}^{(3)}_{\langle d,4m\rangle,\langle1,1+jm\rangle}$ does not depend on $j$.
Therefore from (\ref{R13'}) we get
\begin{equation}\label{R13'est}
R^\prime_{1,3}(X)=0\,.
\end{equation}

From (\ref{R13}), (\ref{R13*est}) and (\ref{R13'est}) it follows that
\begin{equation}\label{R13est}
R_{1,3}(X)\ll X^{2}\mathcal{L}^{-1}\,.
\end{equation}
\textbf{The estimation of} $\textbf{R}_{\textbf{2},\textbf{2}}\textbf{(X).}$
Let $\mathcal{P}$ be the set of primes $X/2<p\leq X$ such that $p-1$ has a divisor lying between $\sqrt{X}\mathcal{L}^{-\omega}<d<\sqrt{X}\mathcal{L}^{\omega}$, (with $\omega=B(10)+C(10)+1$).
Using (\ref{D}), (\ref{r1r2r3}) and (\ref{Rij}) and the inequality $uv\leq u^2+v^2$ we obtain
\begin{align*}
&R_{2,2}(X)\ll\\
&\ll\mathcal{L}^3\sum\limits_{X/2<p_1,p_2,p_3\leq X\atop{p_1+p_3=2p_2\atop{p_1\in\mathcal{P}}}}
\bigg|\sum\limits_{d\mid p_2-1\atop{D<d<X/D}}\chi(d)\bigg|^2
+\mathcal{L}^3\sum\limits_{X/2<p_1,p_2,p_3\leq X\atop{p_1+p_3=2p_2\atop{p_2\in\mathcal{P}}}}
\bigg|\sum\limits_{d\mid p_1-1\atop{D<d<X/D}}\chi(d)\bigg|^2\\
&=\mathcal{L}^3\sum\limits_{X/2<p_2\leq X}\bigg|\sum\limits_{d\mid p_2-1\atop{D<d<X/D}}\chi(d)\bigg|^2
\sum\limits_{p_1+p_3=2p_2\atop{p_1\in\mathcal{P}}}1
+\mathcal{L}^3\sum\limits_{X/2<p_1\leq X}\bigg|\sum\limits_{d\mid p_1-1\atop{D<d<X/D}}\chi(d)\bigg|^2
\sum\limits_{2p_2-p_3=p_1\atop{p_2\in\mathcal{P}}}1\,.
\end{align*}
In order to estimate the sums over $p_1,p_3$  and $p_2,p_3$  we apply respectively
Lemma \ref{Tolev} and Lemma \ref{Tol-Dim} and we get
\begin{equation*}
R_{2,2}(X)\ll X\mathcal{L}^{1-2\theta_0}(\log\mathcal{L})^6
\sum\limits_{X/2<p\leq X}\bigg|\sum\limits_{d\mid p-1\atop{D<d<X/D}}\chi(d)\bigg|^2\,.
\end{equation*}

Using Lemma \ref{Hooley2} we find
\begin{equation}\label{R22est}
R_{2,2}(X)\ll X^2\mathcal{L}^{-2\theta_0}(\log\mathcal{L})^{13}\,.
\end{equation}
\textbf{The estimation of} $\textbf{R}_{\textbf{2},\textbf{3}}\textbf{(X).}$
From (\ref{I1}), (\ref{r1r2r3}) and (\ref{Rij}) we have
\begin{align*}
R_{2,3}(X)&=\sum\limits_{X/2<p_1,p_2,p_3\leq X\atop{p_1+p_3=2p_2}}r_2(p_1-1)\log p_1\log p_2\log p_3
\sum\limits_{m\mid p_2-1\atop{\frac{p_2-1}{m}\geq X/D}}
\chi\left(\frac{p_2-1}{m}\right)\\
&=\sum\limits_{X/2<p\leq X}r_2(p-1)\log p\sum\limits_{m<D\atop{2\mid m}}\sum\limits_{j=\pm1}\chi(j)
I^{(1)}_{4m,1+jm}(p,X,J_m)\,,
\end{align*}
where $J_m=\big(\max\{1+mX/D,X/2\},X\big]$. Using (\ref{Delta1}) we obtain
\begin{equation}\label{R23}
R_{2,3}(X)=R^\prime_{2,3}(X)+R^*_{2,3}(X)\,,
\end{equation}
where
\begin{align}\label{R23'}
&R^\prime_{2,3}(X)=\sum\limits_{X/2<p\leq X}r_2(p-1)\log p\sum\limits_{m<D\atop{2\mid m}}\frac{\Phi^{(1)}(p,X,J_m)}{\varphi(4m)}\sum\limits_{j=\pm1}\chi(j)
\mathfrak{S}_{4m,1+jm}(p)\,,\\
\label{R23*}
&R^*_{2,3}(X)=\sum\limits_{X/2<p\leq X}r_2(p-1)\log p\sum\limits_{m<D\atop{2\mid m}}
\sum\limits_{j=\pm1}\chi(j)\Delta^{(1)}_{4m,1+jm}(p,X,J_m)\,.
\end{align}

Consider $R^\prime_{2,3}(X)$. From Lemma \ref{j} we have that $\mathfrak{S}_{4m,1+jm}(p)$ does not depend on $j$.
Therefore using (\ref{R23'}) we find
\begin{equation}\label{R23'est}
R'_{2,3}(X)=0\,.
\end{equation}

Next we consider $R^*_{2,3}(X)$. From (\ref{r1r2r3}), (\ref{R23*}) and Cauchy's inequality we get
\begin{align}\label{R23*est1}
R^*_{2,3}(X)&\ll\mathcal{L}\sum\limits_{X/2<p\leq X}\tau(p-1)\sum\limits_{m<D\atop{2\mid m}}
\sum\limits_{j=\pm1}\big|\Delta^{(1)}_{4m,1+jm}(p,X,J_m)\big|\nonumber\\
&\ll\mathcal{L}\sum\limits_{X/2<n\leq X}\tau(n)\sum\limits_{m<D\atop{2\mid m}}
\sum\limits_{j=\pm1}\big|\Delta^{(1)}_{4m,1+jm}(n,X,J_m)\big|\nonumber\\
&\ll\mathcal{L}U_1^{1/2}V_1^{1/2}\,,
\end{align}
where
\begin{align*}
&U_1=\mathcal{L}\sum\limits_{X/2<n\leq X}\tau^2(n)\sum\limits_{m<D\atop{2\mid m}}
\sum\limits_{j=\pm1}\big|\Delta^{(1)}_{4m,1+jm}(n,X,J_m)\big|\,,\\
&V_1=\mathcal{L}\sum\limits_{X/2<n\leq X}\sum\limits_{m<D\atop{2\mid m}}
\sum\limits_{j=\pm1}\big|\Delta^{(1)}_{4m,1+jm}(n,X,J_m)\big|\,.
\end{align*}
We use the trivial estimate $\Delta^{(1)}\ll\mathcal{L}^2Xm^{-1}$ and the inequality
$\sum_{n\leq y}\tau^2(n)\ll y\log^3y$ to obtain
\begin{equation}\label{U1est}
U_1\ll X^2\mathcal{L}^6\,.
\end{equation}
We estimate $V_1$ using (\ref{D}) and Lemma \ref{Laporta} and we find
\begin{equation}\label{V1est}
V_1\ll X^2\mathcal{L}^{-10}\,.
\end{equation}
From (\ref{R23*est1}) -- (\ref{V1est})  it follows that
\begin{equation}\label{R23*est2}
R^*_{2,3}(X)\ll X^2\mathcal{L}^{-1}\,.
\end{equation}

Now bearing in mind (\ref{R23}), (\ref{R23'est}) and (\ref{R23*est2}) we obtain
\begin{equation}\label{R23est}
R_{2,3}(X)\ll X^2\mathcal{L}^{-1}\,.
\end{equation}
\textbf{The estimation of} $\textbf{R}_{\textbf{3},\textbf{3}}\textbf{(X).}$
From (\ref{I3}), (\ref{r1r2r3}) and (\ref{Rij}) we have
\begin{align*}
R_{3,3}(X)&=\sum\limits_{X/2<p_1,p_2,p_3\leq X\atop{p_1+p_3=2p_2}}\log p_1\log p_2\log p_3
\sum\limits_{m_1\mid p_1-1\atop{\frac{p_1-1}{m_1}\geq X/D}}
\chi\left(\frac{p_1-1}{m_1}\right)
\sum\limits_{m_2\mid p_2-1\atop{\frac{p_2-1}{m_2}\geq X/D}}
\chi\left(\frac{p_2-1}{m_2}\right)\\
&=\sum\limits_{m_1,m_2<D\atop{2\mid m_1,2\mid m_2}}\sum\limits_{j_1=\pm1\atop{j_2=\pm1}}\chi(j_1)\chi(j_2)
I^{(3)}_{\langle 4m_1,4m_2\rangle,\langle1+j_1m_1,1+j_2m_2\rangle}(X,\textbf{J}_\textbf{m})\,,
\end{align*}
where $\textbf{J}_\textbf{m}=\langle J_{m_1},J_{m_2}\rangle$;\;
$J_{m_\nu}=\big(\max\{1+m_\nu X/D,X/2\},X\big]$,\,$\nu=1,2$.\\
We write
\begin{equation}\label{R33}
R_{3,3}(X)=R^\prime_{3,3}(X)+R^*_{3,3}(X)\,,
\end{equation}
where
\begin{align}\label{R33'}
&R^\prime_{3,3}(X)=\sum\limits_{m_1,m_2<D\atop{2\mid m_1,2\mid m_2}}\frac{\Phi^{(3)}(X,\textbf{J}_\textbf{m})}{\varphi(4m_1)\varphi(4m_2)}\sum\limits_{j_1=\pm1\atop{j_2=\pm1}}\chi(j_1)\chi(j_2)
\mathfrak{S}^{(3)}_{\langle 4m_1,4m_2\rangle,\langle1+j_1m_1,1+j_2m_2\rangle}\,,\\
\label{R33*}
&R^*_{3,3}(X)=\sum\limits_{m_1,m_2<D\atop{2\mid m_1,2\mid m_2}}\sum\limits_{j_1=\pm1\atop{j_2=\pm1}}\chi(j_1)\chi(j_2)
\Delta^{(3)}_{\langle 4m_1,4m_2\rangle,\langle1+j_1m_1,1+j_2m_2\rangle}(X,\textbf{J}_\textbf{m})\,.
\end{align}

Consider first $R^\prime_{3,3}(X)$. According to Lemma \ref{j} the expression
$\mathfrak{S}^{(3)}$ does not depend on $j_2$.
Therefore from (\ref{R33'}) it follows that
\begin{equation}\label{R33'est}
R^\prime_{3,3}(X)=0\,.
\end{equation}

Consider now $R^*_{3,3}(X)$.
Using (\ref{D}) and Lemma \ref{Peneva,Tolev} we find
\begin{equation}\label{R33*est}
R^*_{3,3}(X)\ll X^{2}\mathcal{L}^{-1}\,.
\end{equation}

Now taking into account (\ref{R33}), (\ref{R33'est}) and (\ref{R33*est}) we obtain
\begin{equation}\label{R33est}
R_{3,3}(X)\ll X^{2}\mathcal{L}^{-1}\,.
\end{equation}
\textbf{The estimation of} $\textbf{R}_{\textbf{2},\textbf{1}}\textbf{(X).}$
Using (\ref{I1}), (\ref{Delta1}), (\ref{r1r2r3}) and (\ref{Rij}) we write
\begin{equation*}
R_{2,1}(X)=\sum\limits_{X/2<p\leq X}r_2(p-1)\log p\sum\limits_{d\leq D}\chi(d)I^{(1)}_{d,1}(p,X)=
R^\prime_{2,1}(X)+R^*_{2,1}(X)\,,
\end{equation*}
where
\begin{align*}
&R^\prime_{2,1}(X)=\sum\limits_{X/2<p\leq X}r_2(p-1)\log p\sum\limits_{d\leq D}
\frac{\chi(d)}{\varphi(d)}\mathfrak{S}_{d,1}(p)\Phi^{(1)}(p,X)\,,\\
&R^*_{2,1}(X)=\sum\limits_{X/2<p\leq X}r_2(p-1)\log p\sum\limits_{d\leq D}\chi(d)\Delta^{(1)}_{d,1}(p,X)\,.
\end{align*}

Further arguing as in $R_{1,2}(X)$ we find
\begin{equation}\label{R21est}
R_{2,1}(X)\ll X^2\mathcal{L}^{-\theta_0}(\log\mathcal{L})^7\,.
\end{equation}
\textbf{The estimation of} $\textbf{R}_{\textbf{3},\textbf{2}}\textbf{(X).}$
From (\ref{I2}), (\ref{r1r2r3}) and (\ref{Rij}) we have
\begin{align*}
R_{3,2}(X)&=\sum\limits_{X/2<p_1,p_2,p_3\leq X\atop{p_1+p_3=2p_2}}r_2(p_2-1)\log p_1\log p_2\log p_3
\sum\limits_{m\mid p_1-1\atop{\frac{p_1-1}{m}\geq X/D}}
\chi\left(\frac{p_1-1}{m}\right)\\
&=\sum\limits_{X/2<p\leq X}r_2(p-1)\log p\sum\limits_{m<D\atop{2\mid m}}\sum\limits_{j=\pm1}\chi(j)
I^{(2)}_{4m,1+jm}(p,X,J_m)\,,
\end{align*}
where $J_m=\big(\max\{1+mX/D,X/2\},X\big]$.

Further working as in $R_{2,3}(X)$ we obtain
\begin{equation}\label{R32est}
R_{3,2}(X)\ll X^2\mathcal{L}^{-1}\,.
\end{equation}
\textbf{The estimation of} $\textbf{R}_{\textbf{3},\textbf{1}}\textbf{(X).}$
Arguing similar to $R_{1,3}(X)$ we find
\begin{equation}\label{R31est}
R_{3,1}(X)\ll X^2\mathcal{L}^{-1}\,.
\end{equation}
\textbf{The end of the proof.} The asymptotic formula (\ref{RAs}) follows from (\ref{R16}), (\ref{Rij}),
(\ref{R11est}), (\ref{R12est}), (\ref{R13est}), (\ref{R22est}), (\ref{R23est}), (\ref{R33est}),
(\ref{R21est}), (\ref{R32est}) and (\ref{R31est}).

The Theorem is proved.

\vspace{12pt}
\baselineskip10pt

\vskip10pt

\begin{flushleft}
{\it Received on December, 2015} \\
{\it Revised on December, 2015} \\
\end{flushleft}

\vskip20pt
\footnotesize
\begin{flushleft}
S. I. Dimitrov\\
Faculty of Applied Mathematics and Informatics\\
Technical University of Sofia \\
8, St.Kliment Ohridski Blvd. \\
1756 Sofia, BULGARIA\\
e-mail: sdimitrov@tu-sofia.bg\\
\end{flushleft}
\end{document}